\documentclass[11pt]{article}

\usepackage{amssymb}
\usepackage{amsfonts}
\usepackage{amsmath,amsthm}
\usepackage{latexsym}
\usepackage[numbers,sort&compress]{natbib}
\usepackage{array}
\usepackage{tikz,color,comment}
\usepackage{times}
\usepackage{enumitem}
\usepackage{authblk}
\usepackage{hyperref}

\newcommand{\R}{\mathbb{R}}

\newcommand{\Gnp}{{\mathbb G}_{n,p}}

\newcommand{\pr}{{\mathbb P}}
\newcommand{\E}{{\mathbb E}}

\newcommand{\N}{{\mathbb{N}}}
\newcommand{\Bin}{{\mathrm{Bin}}}

\usepackage[normalem]{ulem}

\newenvironment{enumerate*}%
  {\vspace{-0.3cm}%
  \begin{enumerate}%
    \setlength{\itemsep}{0pt}%
    \setlength{\parskip}{0pt}}%
  {\end{enumerate}}

\newenvironment{itemize*}%
  {\vspace{-0.3cm}%
  \begin{itemize}%
    \setlength{\itemsep}{0pt}%
    \setlength{\parskip}{0pt}}%
  {\end{itemize}}

\newcommand{\mult}{\ensuremath{\mathrm{mult}}}
\newcommand{\vset}{\ensuremath{V}} 
\newcommand{\eset}{\ensuremath{E}}
\newcommand{\forest}{\mathcal{F}}
\newcommand{\I}[1]{\ensuremath{\mathbf{1}_{\{#1\}}}}
\newcommand{\eqdist}{\ensuremath{\stackrel{d}{=}}}
\newcommand{\convdist}{\ensuremath{\stackrel{d}{\rightarrow}}}
\newcommand{\graphs}{\mathcal{G}}

\numberwithin{equation}{section}

\title{\vspace{-10mm}Scaling limits of random graphs}
\author[1]{Louigi Addario-Berry}
\author[2]{Christina Goldschmidt}
\affil[1]{Department of Mathematics and Statistics, McGill University. \\louigi@gmail.com}
\affil[2]{Department of Statistics, University of Oxford. \\goldschm@stats.ox.ac.uk}

\begin{document}

\maketitle

\begin{abstract}
This work will appear as a chapter in a forthcoming volume titled ``Topics in Probabilistic Graph Theory'.
A theory of scaling limits for random graphs has been developed in recent years. This theory gives access to the large-scale geometric structure of these random objects in the limit as their size goes to infinity, with distances appropriately rescaled. We start with the simplest setting of random trees, before turning to various examples of random graphs, including the critical Erd\H{o}s--R\'enyi random graph. 
\end{abstract}

\renewcommand{\baselinestretch}{0.5}\normalsize
\tableofcontents
\renewcommand{\baselinestretch}{1.00}\normalsize

\section{Introduction}

\setcounter{section}{1}
\setcounter{equation}{0}

\noindent A sequence of random variables $(Z_n)_{n \ge 1}$ has a \emph{scaling limit} if there exists a sequence $(a_n)_{n \ge 1}$ of real numbers such that $(Z_n/a_n)_{n \ge 1}$ converges in distribution to some limit random variable $Z$.  Perhaps the most famous example of a scaling limit is 

\medskip

\noindent \textbf{Theorem 1.1} (The central limit theorem) \
 \emph{Suppose that $X_1, X_2, \ldots$ are independent and identically distributed random variables with mean $0$ and finite variance $\sigma^2$.  If $Z_n = X_1 + X_2 + \cdots + X_n$ and $a_n = \sigma \sqrt{n}$, then $Z_n/a_n \overset{d}{\to} Z$ where $Z \sim N(0,1)$.}

\medskip

Roughly speaking, we may say that, for large $n$, the distribution of $Z_n/a_n$ is approximately $N(0,1)$.  But if, rather than having finite variance, the random variables $(X_i)_{i \ge 1}$ are such that $\mathbb{P}(X_i \ge x)=(1+o(1))c/x^{\alpha}$ and $\mathbb{P}(X_i \le- x)=(1+o(1)) c'/x^{\alpha}$ as $x \to \infty$, for some $\alpha \in (1,2)$ and $c,c' \in \R$, then $Z_n/n^{1/\alpha}$ has a distributional limit which is a so-called {\em $\alpha$-stable distribution} and is uniquely determined by $\alpha$, $c$ and $c'$. If the distribution of the summands of $Z_n$ is allowed to depend on $n$, then even more scaling limits are possible; the collection of distributions which can arise in this way are the {\em infinitely divisible distributions}.

Scaling limits are ubiquitous in modern probability, not only for sequences of real-valued random variables, but also for  random objects with values in more general spaces. For a more complex (but still classical) example, define $B^{(n)}(t)=Z_{nt}/{a_n}$ for $t\in \N/n$, and extend this to general $t \ge 0$ by linear interpolation to define a continuous process $B^{(n)}=(B^{(n)}(t),t \ge 0)$. If the summands $(X_i)_{i \ge 1}$ still have mean 0 and have finite variance, then {\em Donsker's theorem} states that $B^{(n)}$ converges in distribution to a standard Brownian motion. If the summands do not have finite variance, or are allowed to depend on $n$, then other limits are possible, and the collection of processes that can arise in this way are called {\em L\'evy processes}.

The theory of scaling limits of random trees and random graphs has developed over the last 35 years or so, but encompasses (and unifies) research that dates back much further. We aim to take a sequence of random trees or graphs of increasing size $n$ (typically measured in terms of the number of vertices) and describe their large-scale metric behaviour, after a suitable rescaling, as $n \to \infty$.  In particular, it is important to identify the correct length scaling for the random graph in question.  This perspective lends itself particularly well to the study of global geometric properties, such as typical distances, diameters, etc., in sparse graph models. On the other hand, it does not capture local properties, such as degrees or subgraph counts.  Indeed, the limiting objects we  discuss are not themselves trees or graphs in the usual combinatorial senses of those words, but are rather continuum analogues in which we can, for example, no longer distinguish a vertex and its adjacent edges.  The right mathematical framework for the limit objects turns out to be compact metric spaces, often endowed with a probability measure which provides a mechanism for sampling points.  The limit objects we uncover are fascinating in their own right; for example, several are random fractals with particularly nice distributional properties.  We do not attempt a comprehensive overview of what is by now an extensive literature, but rather aim to give a flavour of the results and the mathematical tools involved in their proofs.

We begin this chapter with an account of the scaling limit theory for random trees, which started with pioneering work of Aldous, Le Gall, and others on the Brownian continuum random tree, which is the scaling limit of most `uniform-like' random trees, such as uniform random labelled trees and uniform random unlabelled trees.  In each of these cases, it turns out that if the tree has $n$ vertices, then the typical length-scaling is $\sqrt{n}$. The fact that the same scaling limit arises for different discrete models is an important feature of the theory, and is usually referred to as \emph{universality}, a term borrowed from physics. (The central limit theorem is a universal result in that the precise details of the distribution of the underlying i.i.d.\ random variables $X_1, X_2, \ldots$ are unimportant, as long as they are centred and have finite variance.) We also discuss briefly some generalizations to non-uniform models of trees, including those with heavy-tailed degree distributions, for which different scaling limits can arise. 

We then turn to random graphs.  Here, the scaling limit perspective turns out to be particularly valuable in the context of random graphs at the point of the phase transition for the emergence of a giant component.  We will focus on the Erd\H{o}s--R\'enyi random graph $\Gnp$, although the phenomena we describe turn out to be rather generic for the most commonly studied models of random graphs.  

Suppose, then, that we take $\Gnp$ with $p = c/n$, for some constant $c > 0$ which may be thought of as approximately the average degree.  If $c < 1$, then the largest component of $\Gnp$ has size $O_p(\log n)$, whereas for $c > 1$ there exists a unique component of size $\Theta_p(n)$ (and all of the others are again $O_p(\log n)$). Indeed, for any $c \ne 1$, the size of the largest component is concentrated, in the sense that when divided by $n$ it converges in probability to a constant. The small components, both in the sub- and super-critical regimes are essentially trees, whereas the giant component is a much denser object, containing cycles with lengths of the order of $\log n$, and does not turn out to yield a meaningful scaling limit.  

In the critical phase where $c = 1$ (or, indeed, inside the \emph{critical window}, for which $p = 1/n + \lambda n^{-4/3}$, for $\lambda \in \R$), the largest components have sizes of the order of $n^{2/3}$, and there are a diverging number of them.  On rescaling, the component sizes retain some randomness (see Theorem 9.1 below) and are tree-like, in the sense that each component has finitely many edges more than a tree. Moreover, typical distances between vertices are of the order of the square root of the size, as for uniform trees. These largest critical components turn out to have interesting scaling limits, which we describe in various different ways, including via a classical core/kernel decomposition.

We observe at this point that several other notions of graph limits have been developed in recent years.  The graphs we consider in this chapter, which all have bounded average degrees, may also be fruitfully studied via \emph{local weak convergence} (or \emph{Benjamini--Schramm convergence}) in which, rather than rescaling, we instead focus on the ball of radius $R$ in the graph distance around a particular point (often chosen uniformly at random); if this ball converges in distribution for every $R > 0$ then we obtain a \emph{local weak limit} which is typically an infinite (but locally finite) graph. In some sense, this perspective is complementary to ours, in that it captures precisely the sort of local subgraph behaviour that gets washed away in the scaling limit.  Other notions are better adapted to the setting of dense graphs -- for example, the theory of graphons.

A recurring theme in this chapter is that our ability to analyse a particular model of random trees or graphs is intimately related to having a good algorithm for generating it.  To take an example, it is straightforward to make sense of a uniform random labelled tree on $n$ vertices: Cayley's formula tells us that there are $n^{n-2}$ labelled trees on $n$ vertices, and we simply pick one of them uniformly at random, but this tells us nothing about how we might reasonably \emph{generate} such a tree. In this case, there is a beautiful algorithm which provides a way to generate the tree vertex by vertex, and which also hints at a limiting version which actually builds the Brownian continuum random tree.

\section{Models of random trees and graphs}

\setcounter{section}{2}
\setcounter{equation}{0}

\noindent Perhaps the most probabilistically natural model of a random tree describes the genealogy of a branching process.  (These are sometimes called `Galton--Watson trees' in the literature, although we prefer the term \emph{Bienaym\'e tree}.)  We fix an offspring distribution with probability mass function $(p_k)_{k \ge 0}$, and start from a population consisting of a single individual. Each individual in the population independently gives birth to a number of children distributed as $(p_k)_{k \ge 0}$.  If it ever happens that there are no individuals left, then we say that the population has become extinct.  Let us suppose that the offspring distribution is non-trivial, in the sense that $p_0 + p_1 < 1$, and let $m = \sum_{k \ge 1} k p_k$ be the mean of the offspring distribution. It is then a standard result that the population becomes extinct with probability 1 if $m \le 1$, and has positive probability of non-extinction if $m > 1$.

The genealogy of such a branching process is most naturally formulated as a tree which is rooted at the original individual, and is ordered so that we have a first child of the root, a second child of the root, and so on.  In this chapter, we are mostly interested in trees of a fixed size, and so we will typically consider a branching process tree conditioned to have $n$ vertices.  (We usually assume, for simplicity of presentation, that this is an event of positive probability. This will not be the case for every offspring distribution -- consider, for example, the situation where $p_0+p_2 = 1$ -- but the only adjustment necessary is to consider a subsequence of integers for which the event does have positive probability.)  In the cases where $m < 1$ and $m > 1$, the probability that the tree has precisely $n$ vertices is exponentially small in $n$, whereas if $m = 1$ it often has probability which decays roughly like an inverse power of $n$ (with the power depending on the tail of the offspring distribution). So we will usually be interested in the critical case.

This gives a very rich class of random trees to deal with.  At first sight, it may not look natural from the combinatorial perspective. However, it turns out that in many cases, trees which are chosen uniformly at random from a given combinatorial class may be cast into this framework.  For example:

\begin{itemize}[leftmargin=10pt]
\item A uniform random rooted ordered tree with $n$ vertices arises when we take the offspring distribution to be Geometric$(\tfrac{1}{2})$, so that $p_k = \left(\frac{1}{2}\right)^{k+1}$ for $k \ge 0$.
\item A uniform binary rooted ordered tree with $n$ vertices arises when we take the offspring distribution to be $p_0 = \tfrac{1}{2}$, $p_2 = \tfrac{1}{2}$, provided that $n$ is odd.
\item A uniform random rooted labelled tree with $n$ vertices arises when we take the offspring distribution to be Poisson$(1)$, so that $p_k = e^{-1}/k!$ for $k \ge 0$; assign the vertices a uniformly random labelling by $\{1,2, \ldots, n\}$, and then keep the root and forget the ordering.
\end{itemize}

We will also be interested in the setting of rooted labelled trees, where we fix in advance the \emph{child sequence} (that is, the numbers $c_1, c_2, \ldots, c_n$ of children of the vertices $1, 2, \ldots, n$, respectively) and generate a uniformly random rooted labelled tree with that child sequence. (Note that we must have $\sum_{i=1}^n c_i = n-1$, because every vertex is the child of some other vertex except for the root.)  In the next section, we provide a simple algorithm to generate such a tree. We may fit the model of conditioned Bienaym\'e trees into this framework by taking the child sequence to consist of independent and identically distributed random variables (usually with mean 1, so that we are in the critical case) with sum conditioned to equal $n-1$.

Let us turn now to random graphs.  We will typically think of graphs with vertices labelled by the integers $1, 2, \ldots, n$.  The analogue of fixing the child sequence of a tree is to fix the \emph{degree sequence}, $d_1, d_2, \ldots, d_n$ of the vertices $1, 2, \ldots, n$.  In order for this to make sense, we require $\sum_{i=1}^n d_i$ to be even and, more strongly, for the degree sequence to be \emph{graphical} (with at least one simple graph having that degree sequence; explicit conditions for this are given by the Erd\H{o}s--Gallai theorem~\cite{MR144332}).  Here, there is a standard method of generating a \emph{multigraph} with the given degrees, called the \emph{configuration model}. Conditioned on simplicity, this is a uniform random graph with a given degree sequence.  It turns out to be helpful for the analysis if there is an extra layer of randomness, so that we take the degrees themselves to be independent and identically distributed random variables.  

Perhaps the simplest sparse random graph, the Erd\H{o}s--R\'enyi random graph $\Gnp$ for $p = c/n$, may (roughly speaking) be thought of as a uniform random graph with i.i.d.\ Poisson$(c)$ random degrees.  However, the extra independence between different edges inherent in the usual description of $\Gnp$ makes it particularly tractable.  We observe that a component of $\Gnp$ with vertex-set $\{v_1, v_2, \ldots, v_k\}$ and $m$ edges is a uniform random connected graph with $m$ edges; in particular, if $m=k-1$, it is a uniform random labelled tree with vertex-set $\{v_1, v_2, \ldots, v_k\}$.

\section{Connected graphs as metric spaces, and convergence}

\setcounter{section}{3}
\setcounter{equation}{0}

\noindent As mentioned above, a convenient formulation of our scaling limit results requires us to think of a connected graph as a metric space.  This is quite natural: the points of the metric space are the vertices of the graph, and the metric is given by the graph distance -- that is, the number of edges in a shortest path between two vertices. We also endow our metric spaces with a probability measure, which we always take to be the uniform measure on the vertices.

In order to talk about convergence, we need a topology, and this is provided by the so-called Gromov--Hausdorff--Prokhorov (GHP) distance on the set of measured metric spaces -- that is, metric spaces equipped with a measure -- up to measure-preserving isometry.  The details are quite technical, and we omit them here (see \cite{MR1835418} for a slightly simpler concept, the Gromov--Hausdorff distance, which ignores the measured structure, and \cite{MR3522292} and \cite[Section {6}]{MR4712854} for the full version).  Instead, we give a pair of sufficient conditions for convergence which are easier to understand.

The basic philosophy is that we can `find' the whole of a well-behaved metric space by repeatedly sampling uniformly at random from it. In order to reconstruct the metric space, it is sufficient to know all the pairwise distances between those uniform vertices.  If the matrix of pairwise distances converges, and we have an appropriate tightness condition, then we get convergence of the measured metric spaces.

Suppose that, for each $n \in \N$, $(X^n, d^n)$ is a random compact metric space and $\mu^n$ is a Borel probability measure on $X^n$.  For each $n \in \N$, conditionally on $(X^n,d^n,\mu^n)$, let $(u_i^n)_{i \ge 1}$ be a sequence of i.i.d.\ samples drawn from $X^n$, according to the distribution $\mu^n$.  Let $\mathbf{M}^n = (d^n(u_i^n, u^j_n))_{i,j \in \N}$ be the matrix of pairwise distances between the random points of $(X^n, d^n, \mu^n)$.  We say that the sequence $(X^n,d^n,\mu^n)_{n \ge 1}$ is \emph{tight} if, for all $\delta > 0$, we have
\[
\lim_{k \to \infty} \limsup_{n \to \infty} \pr \left(d_H(X^n, \{u_i^n\}_{1 \le i \le k}) > \delta \right) = 0,
\]
where $d_H$ is the usual Hausdorff distance between compact subsets of a metric space.  
(Note that this is stronger than requiring $\mu^n$ to have support $X^n$ for all $n$.) If $\mathbf{M}^n \convdist \mathbf{M}$ as $n \to \infty$ and the sequence is tight, then
\[
(X^n, d^n, \mu^n) \convdist (X, d, \mu) \quad \text{in the GHP sense,}
\]
 for some measured metric space $(X,d,\mu)$ which is such that, if $(u_i)_{i \ge 1}$ are i.i.d.\ samples from $\mu$, then $(d(u_i, u_j))_{i,j \in \N} \eqdist \mathbf{M}$ and $X = \mathrm{supp}(\mu)$.

\section{The contour process of a tree}

\setcounter{section}{4}
\setcounter{equation}{0}

\noindent In this section, it is useful to treat finite trees $T$ as metric spaces by viewing each of their edges as an isometric copy of a unit interval. 

Given $\zeta > 0$, an {\em excursion} of length $\zeta$ is a continuous function $f:[0,\zeta]\to [0,\infty]$ with $f(0)=f(\zeta)=0$. 
To a finite rooted ordered tree $T$ with $n$ vertices, we associate an excursion $e=e_T$ of length $2n-2$, called the {\em contour process} of $T$, as follows. Starting from the root, explore the contour of the tree in clockwise order at unit speed, finishing the first time every edge has been traversed twice (once in each direction). Then set $e_T(x)$ to be the distance to the root at time $x$ of the exploration; see Fig.\ 1(a) and 1(b) for an example.

\begin{center}
\begin{tabular}{ccc}
\includegraphics[width=4cm,page=1]{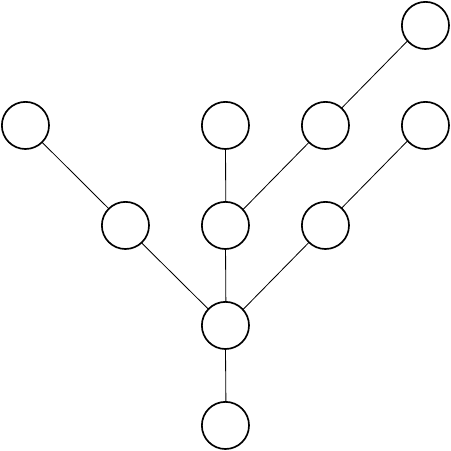} & \hspace{1cm} & 
\includegraphics[width=6cm,page=2]{contourexamplesplitup.pdf} \\
(a) & \hspace{1cm} & (b) \\
\end{tabular}
\smallskip

\noindent {\small Fig.\ 1. (a) A rooted ordered tree $T$ with $10$ vertices; (b) the contour process of $T$}
\end{center}

The metric structure of the tree $T$ can be recovered from $e_T$ as follows. Given {\em any} excursion $e$ of length $\zeta>0$, first
define a pseudometric on $[0,\zeta]$ by setting 
\[
d_e^0(x,y)=e(x)+e(y)-2\min_{x\wedge y \le z \le x \vee y} e(z),
\]
for $0 \le x,y \le \zeta$. Then define an equivalence relation $\sim$ on $[0,\zeta]$ by setting $x \sim y$ if and only if $d_e^0(x,y)=0$. Finally, let $d_e$ be the metric induced on $[0,\zeta]/\!\sim$ from  $d_e^0$. It is not hard to verify that the result is a compact metric space, $T_{e}=([0,\zeta]/\!\sim,d_e)$. Its construction can be understood visually as follows: imagine applying glue to the underside of the excursion $e$, and then squashing the excursion together horizontally; any two points that meet stick together (which means they are identified under $\sim$), and this gives a tree. 

We also observe that there is a natural way to put a `uniform' probability measure $\mu_e$ on $T_e$,  by simply taking the push-forward of the uniform probability measure on $[0,\zeta]$ onto the quotient space.  It is straightforward to see that if $e=e_T$ is the contour process of a finite tree $T$, then $T_e$ is isometric to $T$ (where we again view the edges of $T$ as unit-length line-segments), and $\mu_e$ is indeed the uniform measure on $T_e$. In fact, the root and the ordered structure can also be recovered, but we do not need this.

The connection between trees and excursions has played a fundamental role in the study of scaling limits of random trees. It is perhaps most clearly exhibited by the case of a uniform random rooted ordered tree of size $n$. Recall that such a tree $T_n$ is equivalent to a Geometric$(\tfrac{1}{2})$ Bienaym\'e tree, conditioned to have size $n$. In this case, the contour exploration $e_n=e_{T_n}$ is precisely a symmetric simple random walk path of length $2n-2$, conditioned to stay non-negative and to return to $0$ at time $2n-2$. The scaling limit of such a conditioned path is a {\em Brownian excursion} $\mathbf{e}$, which is a Brownian motion path run for time $1$, conditioned to stay non-negative and to return to the origin at time $1$. (This is a degenerate conditioning, but it can be made rigorous.) The tree $\mathcal{T}=T_{\mathbf{e}}$, built from a Brownian excursion, is the {\em Brownian continuum random tree} (CRT) mentioned in the introduction. It turns out that the contour process $e_n$ of any conditioned Bienaym\'e tree with critical finite variance offspring distribution also has $\mathbf{e}$ as its scaling limit, and so we say that such offspring distributions are in the {\em Brownian domain of attraction}. (In general, the contour process of a Bienaym\'e tree does not have a nice distribution, but below we sketch the argument that gives the convergence.)  A consequence is that $\mathcal{T}$ is the scaling limit for all of these trees. The following theorem is due to Aldous~\cite{MR1207226} and Le Gall~\cite{MR2203728}.

\medskip
\noindent \textbf{Theorem 4.1} \
\emph{Let $T_n$ be a conditioned Bienaym\'e tree with critical  offspring distribution of finite variance $\sigma^2 \in (0,\infty)$. Let $d_n$ be the graph distance on $T_n$ and $\mu_n$ be the uniform measure on its vertices. Then
\[
(V(T_n), \tfrac{1}{2} \sigma n^{-1/2}d_n, \mu_n) \convdist (T_{\mathbf{e}}, d_{\mathbf{e}}, \mu_{\mathbf{e}}),
\]
in the GHP sense, as $n \to \infty$.}

\medskip

The Brownian CRT is also the scaling limit for numerous combinatorial models which cannot be described via branching processes, including some non-tree models such as random outerplanar graphs, which are `globally tree-like' in spite of possessing cycles; the cycles disappear `in the scaling limit'.

We conclude this section with an observation about the excursions of a Brownian motion which we shall need later -- namely, that they possess a \emph{scaling property}. Indeed, if we write $\mathbf{e}_x$ for a Brownian excursion of length $x$, then we have
\begin{equation}\label{eq:bscale}
x^{-1/2} (\mathbf{e}_x(xt), 0 \le t \le 1) \eqdist (\mathbf{e}(t), 0 \le t \le 1).
\end{equation}
It follows that a Brownian excursion of length $x$ encodes a tree which is distributed as $\mathcal{T}$ with all its distances multiplied by $\sqrt{x}$.

\section{Line-breaking constructions of trees}

\setcounter{section}{5}
\setcounter{equation}{0}

\noindent A more combinatorial approach to the construction of random trees and their limits is hinted at by the following proof of Cayley's formula. The proof exhibits a bijection between the set of labelled rooted trees with vertex-set $[n]$ and the set $[n]^{n-1}$, defined by
\[
\{v_1 v_2 \ldots v_{n-1}: v_i \in [n], 1 \le i \le n-1\}.
\] 

Given such a tree $T$, we first construct an increasing sequence $(S_i(T))_{i \ge 0}$ of subtrees of $T$. The first subtree $S_0(T)$ consists of the root alone. Given the subtree $S_i=S_i(T)$, let $y_{i+1}$ be the lowest-numbered vertex of $T$ which does not belong to $S_i$, let $P_{i+1}$ be the path from $S_i$ to $y_{i+1}$, and form $S_{i+1}$ by adding $P_{i+1}$ to $S_i$. We conclude at the step $I$ when $S_I=T$, and then form an element of $[n]^{n-1}$ as follows: let $P_i^-$ be $P_i$ with its final element removed, so if $P_i=p_1 p_2 \ldots p_l$, then $P_i^-=p_1 p_2 \ldots p_{l-1}$. The coding word for $T$ is then $v_1 v_2 \ldots v_{n-1} = P_1^- P_2^-\ldots P_I^-$.
The length of the coding word is indeed $n-1$, since the number of times that an integer $v$ appears is precisely the number of children of $v$ in $T$. 
An example appears in Fig.\ 2.

\begin{center}
\includegraphics[width=0.4\textwidth]{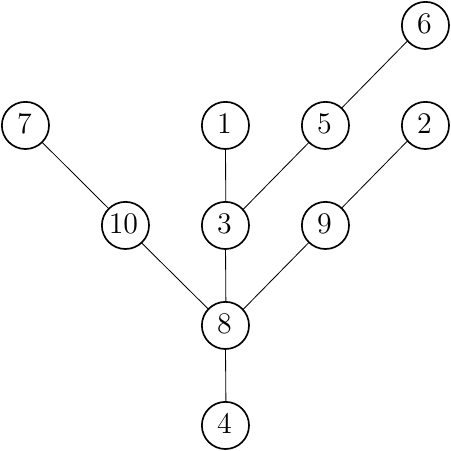}
\end{center}
{\small Fig.\ 2. An example of the line-breaking proof of Cayley's formula. Here, $P_1=4\ 8\ 3\ 1$, $P_2=8\ 9\ 2$, $P_3=3\ 5$, $P_4=5\ 6$, $P_5=8\ 10\ 7$, and $y_1,y_2, \ldots,y_5=1,2,5,6,7$. The coding word for this tree is $4\ 8\ 3\ 8\ 9\ 3\ 5\ 8\ 10$}

\medskip

This construction is reversible. First, if the first entry of the coding word is not~$1$, then $y_1 = 1$; otherwise, $y_1 =2$ (and $T$ has root $1$).
Next, for each $i \ge 1$, we can see where $P_i^-$ ends and $P_{i+1}^-$ begins, because the first element of $P_{i+1}^-$ is either the first repetition in the coding word after $P_i^-$ starts, or  is the element $y_i$. (The identity of $y_i$ can also be deduced, because it is the smallest number that does not appear in $P_1, P_2, \ldots,P_{i-1}$). Moreover, if we restrict to trees rooted at vertex $1$, and if in the Pr\"ufer encoding we use the rule that we always remove the largest-labelled leaf,  then the reversed word $v_{n-1} v_{n-2} \ldots v_2 v_1$ is the Pr\"ufer code of the tree.

This bijection is useful for understanding the distributional properties of a {\em uniformly random rooted tree} $T_n$ with vertex-set $[n]$.  For example, for any $1 \le i <  n$, given the tree $S_i=S_i(T_n)$, 
if $S_i\ne T_n$, then the attachment vertex of the path $P_{i+1}$ to $S_i$ is uniformly distributed over $S_i$. Also, given that $i+1 \not \in S_i$, the length of the path from $i+1$ to $S_i$ is distributed as 
the first time to see either an element of $S_i$, the entry $(i+1)$, or  a repetition, in a sequence of i.i.d.\ Uniform$[n]$ random variables. 
Thus, if $\rho=\rho(T_n)$ denotes the root of $T_n$, then 
\[
\mathbb{P}(\mathrm{dist}_{T_n}(\rho,1)=d) = \frac{d+1}{n}\cdot \prod_{j=1}^{d} \left(1-\frac{j}{n}\right)\, ,
\]
and for $i \ge 1$ and $d \ge 1$,
\[
\mathbb{P}(\mathrm{dist}_{T_n}(i+1,S_i)=d|S_i) = \frac{|S_i|+d}{n}\cdot \prod_{j=1}^{d-1} \left(1-\frac{|S_i|+j}{n}\right)\cdot \mathbf{1}_{\{i+1\not\in S_i\}},
\]
where $\mathbf{1}_E$ is the indicator function of an event or set $E$.
(Note that these are variants of the classical birthday problem, which concerns the time until the first repetition in a sequence of i.i.d.\ Uniform$[n]$ random variables.) One may readily deduce (see \cite{MR1085326}) that, for fixed $i\ge 1$, the subtree $S_i(T_n)$ typically has size $\Theta(n^{1/2})$, and more precisely that
\[
n^{-1/2}(|S_{i}(T_n)|,i \ge 1) \stackrel{\mathop{\mathrm{d}}}{\longrightarrow} (s_i,i \ge 1),
\]
where the limiting vector is the sequence of arrival times of an inhomogeneous Poisson process on $[0,\infty)$ with rate $\lambda(t)=t$. For later reference, we note that, for fixed $k \ge 1$, the joint distribution of $(s_1, s_2, \ldots,s_k)$ is 

\setcounter{equation}{0}
\begin{equation}\label{eq:jointdist}
x_1 x_2 \cdots x_k \cdot e^{-x_k^2/2} \cdot \mathbf{1}_{0 \le x_1 \le x_2 \le \ldots \le x_k}.
\end{equation}

The limiting vector $(s_i, i \ge 1)$ may be used directly to define an increasing sequence of leaf-labelled trees, as follows (see Fig.\ 3). For $i \ge 1$, let $P_i$ be an isometric copy of the interval $[0,s_{i}-s_{i-1}]$, where $s_0=0$. Let $\mathcal{T}_1=P_1$, with the endpoints of the path labelled $0$ and $1$. Then, for $i \ge 1$, form $\mathcal{T}_{i+1}$ from $\mathcal{T}_i$ by identifying one endpoint of $P_{i+1}$ with a point of $\mathcal{T}_i$, distributed according to the (normalized) Lebesgue measure, and giving the other endpoint of $P_{i+1}$  the label $i+1$. The tree $\mathcal{T}_i$ is a binary tree (where all vertices have degrees $1$ or $3$) with root $0$ and non-root leaves labelled $1,2,\ldots,i$. The metric space completion of the limit $\lim_{i \to \infty} \mathcal{T}_i=\bigcup_{i \ge 1} \mathcal{T}_i$ is distributed as the Brownian CRT $\mathcal{T}$. 

\begin{center}
\includegraphics[width=0.75\textwidth]{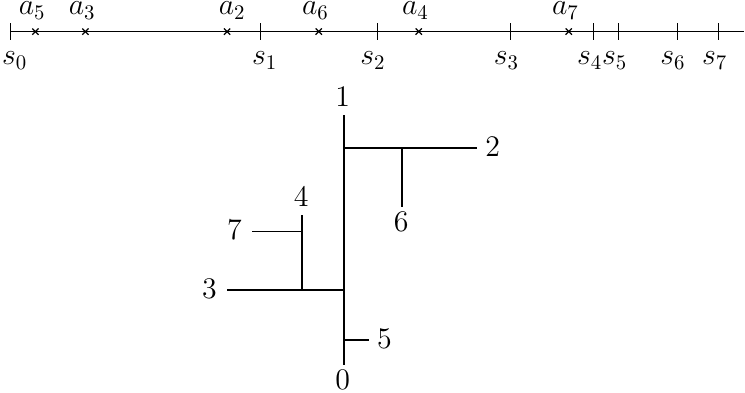}
\end{center}

{\small \noindent Fig.\ 3. An example of the first seven steps of the limiting line-breaking construction. The tree $\mathcal{T}_7$ is depicted below, the points $s_0,s_1,\ldots,s_7$ above are indicated by vertical bars, and for $2 \le i \le 7$, the point $a_i$ of attachment of $P_i$ to $\mathcal{T}_{i-1}$ is indicated by a small $\times$}

\medskip

The above construction of $\mathcal{T}$ is known as Aldous's \emph{line-breaking construction} \cite{MR1085326}. By analogy, we refer to the bijection as the \emph{line-breaking bijection}.  It is in some respects a more versatile construction than those using functional encodings, such as the contour process. In particular, the line-breaking bijection for labelled rooted trees specializes to a bijection between labelled rooted trees in which the number of children of each node is specified and sequences in which the number of times each integer appears is specified. This yields the algorithm promised in Section~2, for generating a random tree with a given child sequence $(c_1,c_2,\ldots,c_n)$ -- simply apply the line-breaking bijection to a uniformly random sequence in which each integer $i\in [n]$ appears $c_i$ times -- and thereby provides access to the distributional properties of a wide range of random tree models. 

Moreover, the limiting line-breaking construction can be directly generalized, in particular to the {\em stable CRTs} first defined by Duquesne and Le Gall \cite{MR1954248}.  (The approach using functional encodings can also be used for the latter purpose, but becomes quite involved.) Stable CRTs are the scaling limits of critical Bienaym\'e trees whose offspring distributions, rather than having finite variance, lie instead in the domain of attraction of a stable random variable. Rather than giving a full definition here, we just observe that this includes any critical offspring distribution, such that $p_k = (1+o(1))c/k^{1+\alpha}$ for $c > 0$ and $\alpha \in (1,2)$.  A line-breaking construction of the stable trees was given in \cite{MR3317158}. 

A remarkable feature of the line-breaking construction of the Brownian CRT is the following. Let $T_i$ be the tree obtained from $\mathcal{T}_i$ by ignoring the edge-lengths and rooting at the leaf with label $0$, and list the edge-lengths of $\mathcal{T}_i$ as $X_{i,1},X_{i,2},\ldots,$ $X_{i,2i-1}$. It then follows fairly straightforwardly from \eqref{eq:jointdist} that the joint density of $T_i$ and $(X_{i,1},X_{i,2},\ldots,X_{i,2i-1})$ is 
\[
f(t;x_1,x_2,\ldots,x_{2i-1}) = s e^{-s^2/2}\, ,
\]
where $s=x_1+x_2+\ldots+x_{2i-1}$ and $t$ is a binary tree with non-root leaves labelled by $1, 2, \ldots, k$ (for a proof, see Aldous~\cite[Lemma 21]{MR1166406}). It follows that the \emph{tree shape} $T_i$ is independent of the edge-lengths, and is uniformly distributed over binary trees with non-root leaves labelled by $1,2,\ldots,k$, and that the edge-lengths are exchangeable. These distributional facts in turn imply that the point of attachment of $P_{i+1}$ to $\mathcal{T}_i$ is uniformly distributed over the edges of $T_i$, so $T_{i+1}$ is built from $T_i$ by subdividing a uniformly random edge and attaching a leaf of label $i+1$ as a child of the new vertex. This rule for sequential generation of leaf-labelled binary trees is known as {\em R\'emy's algorithm} \cite{remy}. 

An extension of R\'emy's algorithm, known as {\em Marchal's algorithm} \cite{MR2508809}, can serve to provide a definition of the stable CRTs mentioned above. Marchal's algorithm is governed by a parameter $\alpha \in (1,2]$, and again generates a nested sequence of trees.

\medskip

\noindent \fbox{
\mbox{\begin{minipage}[t]{0.97\textwidth}
\textsc{Marchal's algorithm}
\begin{itemize}[itemindent=0pt,itemsep=0pt,leftmargin=10pt,topsep=2pt]
\item Start from a tree $T_1^\alpha$ consisting of two adjacent vertices labelled $0$ and $1$, and root the tree at vertex $0$. 
\item For $i =1,2, \ldots$ form $T_{i+1}^{\alpha}$ from $T_i^{\alpha}$ as follows. \\
Assign weight $\alpha - 1$ to the edges of $T_i^\alpha$, and assign weight $c-\alpha$ to any vertex with $c \ge 2$ children. Also, assign weight 0 to vertices with $0$ or $1$ child.

\begin{enumerate*}
\itemindent=-10pt
\item Pick an edge or a vertex with probability proportional to its weight.
\item If a vertex was chosen in step 1, attach a new leaf of label $i+1$ to the chosen vertex. 
\item If an edge was chosen in step 1, subdivide the edge, and attach a new leaf of label $i+1$ to the newly created vertex. 
\end{enumerate*}
\end{itemize}
\end{minipage}}
}

\medskip

The sequence $(T_i^\alpha,i \ge 1)$ converges, after a rescaling of distances by $i^{-(\alpha-1)/\alpha}$, to a limiting random metric space, the {\em $\alpha$-stable CRT}. Moreover, like R\'emy's algorithm, Marchal's algorithm may also be obtained from a line-breaking construction by `ignoring edge-lengths' (see \cite{MR3317158}).

\section{The depth-first queue process of a tree}

\setcounter{section}{6}
\setcounter{equation}{0}

\noindent We conclude our discussion of the scaling limits of random trees with another bijection between trees and lattice paths, which has been fundamental to the study of scaling limits for  both random trees and random graphs. Fix a finite rooted ordered tree $T$, and for $v \in  \vset(T)$, write $c(v)=c(v;T)$ for the number of children of $v$. The \emph{depth-first queue} (DFQ) \emph{process} of $T$ (also called the {\em {\L}ukasiewicz path} of $T$) is defined as follows. List the vertices of $T$ in the order they are first visited by the contour process of $T$ as $v_1, v_2,\ldots,v_n$, and notice that this is a depth-first ordering of $\vset(T)$. Then, for $0 \le i \le n$, set $q_i=q_i(T)=1+\sum_{j=1}^i (c(v_j)-1)$. 
We imagine maintaining a stack (last-in, first-out queue) of vertices to be explored, initially containing just the root. When a vertex $v$ is first visited by the process, it is removed from the stack, and its children (if any) are added. Thus, the net change in stack size is $c(v)-1$.  With this perspective, it is straightforward to verify that at all times  $i$, the quantity $q_i$ gives the total number of {\em younger siblings of ancestors of $v_i$} in $T$. An example appears in Fig.\ 4(a) and 4(b).

\begin{center}
\begin{tabular}{ccc}
\includegraphics[width=4cm,page=1]{contourexamplesplitup.pdf} & \hspace{1cm} & 
\includegraphics[width=6cm,page=3]{contourexamplesplitup.pdf} \\
(a) & \hspace{1cm} & (b) \\
\end{tabular}
\smallskip

{\small \noindent Fig.\ 4. (a) A rooted ordered tree $T$ with $10$ vertices; (b) the DFQ process of $T$}
\end{center}

\medskip

The DFQ process of $T$ is a {\em discrete excursion of length $n$} -- that is, a lattice path $(q_0, q_1,\ldots,q_n)$ with $q_0=1$, $q_n=0$, $q_i > 0$ and $q_{i+1}\ge q_i-1$, for all $0 \le i < n$. Conversely, any discrete excursion of length $n$ is the DFQ process of a unique rooted ordered tree $T$ of size $n$. 
Note that the DFQ process also makes sense for rooted labelled trees, provided that we choose a convention for how to fix a left-to-right order of the children of each vertex.

The DFQ process construction provides a valuable connection between stochastic processes with tractable distributions and random trees and graphs. As a prime example, let $T$ be a critical or subcritical Bienaym\'e tree whose offspring distribution has probability mass function $(p_k)_{k \ge 0}$; then $(q_i(T))$ is a random walk with step distribution $(p_{k+1})_{k \ge -1}$, starting at $1$ and stopped upon hitting $0$. If $T$ is conditioned to have size $n$, then $(q_i(T))$ becomes a random walk excursion of length $n$. In the large-$n$ limit, this provides further connections with continuous stochastic processes; for offspring distributions in the Brownian domain of attraction, such a random walk excursion has a Brownian excursion $\mathbf{e}$ as its scaling limit. On the other hand, if the tree's offspring distribution has a power-law tail, $p_k = (1+o(1))c/k^{1+\alpha}$ for $\alpha \in (1,2)$, then the DFQ process limit is an excursion of a continuous-time process called an \emph{$\alpha$-stable L\'evy process} (whence the name {\em $\alpha$-stable CRT} for the scaling limit of such trees). If the offspring distribution is allowed to depend on~$n$, then other sorts of L\'evy process excursions can also arise. The monograph by Duquesne and Le Gall (\cite{MR1954248}) was foundational in developing this connection between stochastic processes and large random trees.

Since a tree $T$ is reconstructible from either its contour process or its DFQ process, these two process must likewise be constructible one from the other. The details of this are slightly involved, but in certain random cases they become easy to describe in the large-$n$ limit. To give an idea of how it works, we consider the first branch in the line-breaking construction of a random tree $T(c)$ with a given child sequence $c=(c_1, c_2,\ldots,c_n)$. In the coding sequence of such a tree, each integer $i$ appears $c_i$ times, and these appearances occur at uniformly random locations; it follows that the vertices appear in size-biased random order. For example, this means that, for the root of the tree, $\rho=\rho(T(c))$,  
\[
\pr(\rho(T(c))=i) = \frac{c_i}{n-1},
\]
and so 
\[
\mathbb{E}(c_{\rho}) = \frac{1}{n-1}\sum_{i \in [n]} c_i^2.
\]
It is not too hard to see that if $\max_{i \in [n]} c_i = o(n^{1/2})$, then the expected number of children of other vertices along the first branch will also be close to $\frac{1}{n-1}\sum_{i \in [n]} c_i^2$. Moreover, if $vw$ is a parent-child edge of this branch then, by our ordering convention, the left-to-right rank of $w$ is uniform on $\{1,2, \ldots,c_v\}$, and so the number of younger siblings of $w$ is uniform on $\{0,1,\ldots,c_v-1\}$ and thus has mean $\frac{1}{2}(c_v-1)$. Roughly speaking, this yields the result that the expected number of younger siblings of any vertex along the branch is given by
\[
\tfrac{1}{2} \Bigg(\sum_{i \in [n]} \frac{c_i^2}{n-1}-1\Bigg) = \tfrac{1}{2} \sum_{i \in [n]} \frac{c_i(c_i-1)}{n-1}.
\]

For large random trees, the above information is easiest to exploit in the case where $\sum_{i \in [n]} c_i(c_i-1)= \Theta(n)$; in view of the above discussion, this corresponds to a condition that the numbers of children along the first branch are genuinely random (that is, they are not almost all equal to $1$) and have bounded mean. In this case, using the `counting younger siblings of ancestors' characterization of the DFQ process, this information can be converted into a law of large numbers. Writing $v_I$ for the terminal vertex of the first branch in the line-breaking construction, with high probability we then have 
\[
\frac{q_I(T(c))}{\mathrm{dist}(\rho,v_I)} = \tfrac{1}{2}(1+o(1)) \sum_{i \in [n]} \frac{c_i(c_i-1)}{n-1}.
\]
In other words, the heights of the contour process and the DFQ process, when they each visit $v_I$, are related by a roughly constant scaling factor. In fact, it may be shown (see \cite{MR3188597}) that if $\sum_{i\in [n]} c_i(c_i-1)=\Theta(n)$ and $\max_{i \in [n]} c_i = o(n^{1/2})$, then with high probability, 
\[
 q_i(T(c)) =\mathrm{dist}(\rho,v_i) \cdot \tfrac{1}{2} \sum_{i \in [n]} \frac{c_i(c_i-1)}{n-1}+o(n^{1/2}) \ \text{ for all $i \in [n]$.}
\]
Thus, for such trees the two processes are related by an essentially constant scaling factor, related to the empirical variance of the sequence $c$.

Returning to the case of conditioned Bienaym\'e trees, we now have a child sequence consisting of i.i.d.\ random variables conditioned to have sum $n-1$. If the offspring distribution is in the Brownian domain of attraction then certainly the empirical variance of this child sequence is bounded (indeed, it converges as $n \to \infty$ to the variance $\sigma^2$ of the offspring distribution). So we may conclude that the contour process and DFQ process are asymptotically related by a multiplicative factor of $\sigma^2/2$. Since the DFQ process has the Brownian excursion as its scaling limit, so also must the contour process. This completes the promised sketch of the convergence of the contour process in the Brownian case.

In non-Brownian cases (such as that of stable trees), the largest-degree vertices have a macroscopic influence on the tree, and the connection between the DFQ process and the contour process is much looser. Although the transformation relating them is still deterministic (both for finite trees and in continuum trees, for which it is possible to make sense of both processes in many cases), it is no longer simply a matter of applying a constant scaling factor. This is reflected in the fact that, when $\sum_{i \in [n]} c_i^2$ is much larger than $n$, the partial sums of degrees along the first path in the line-breaking construction look asymptotically like a random walk whose step size has infinite mean. For such a walk, the law of large numbers no longer holds, and the value, which is random to first order, is instead dominated by the contributions of the largest summands.

\section{Graphs: cores, kernels and line-breaking constructions}

\setcounter{section}{7}
\setcounter{equation}{0}

\noindent As in the case of trees, there are several possible descriptions of the metric structure of graphs, each with its own strengths. For the moment we focus on connected graphs, because a general graph can be described by specifying the structure of each of its connected components. 

In this section it is useful to allow multigraphs, which may also have loops. For a multigraph $G=(\vset(G),\eset(G))$ and an edge $e \in \eset(G)$, we write $\mult(e)=\mult(e;G)$ for the number of copies  of $e$ in $G$. However, in what follows, the word `graph'  means `simple graph', unless otherwise stated. 

The {\em surplus} of a connected multigraph $G$ is the integer $s(G) = 1+|\eset(G)|-|\vset(G)|$. We write $\graphs_n$ for the set of connected graphs with vertex-set $[n]$, and $\graphs_{n}^s$ for the subset of $\graphs_n$ whose elements have surplus $s$. (Observe that $\mathcal{G}_n^0$ is the set of trees labelled by $[n]$.) It turns out that the global metric structure of a random element of $\graphs_{n}^s$ has a particularly nice asymptotic (large-$n$) description. In order to give this description, we must first introduce a standard decomposition of the graph.

We define a {\em core} to be a connected graph with minimum degree at least $2$. Given a connected graph $G$, the {\em core of $G$}, denoted by $C(G)$, is the maximum induced subgraph of $G$ with minimum degree at least $2$.  Note that $C(G)$ has the same surplus as the original graph $G$. The {\em kernel} of $G$ is the multigraph $K(G)$ obtained from the core $C(G)$ by replacing any path of maximal length whose internal vertices all have degree $2$ by a single edge. Note that $K(G)$ can have multiple edges and loops (see Fig.\ 5). If $G$ has surplus $1$, the kernel is empty by convention.

\begin{center}
\begin{tabular}{ccc}
\parbox[c]{5cm}{\includegraphics[page=1,scale=0.6]{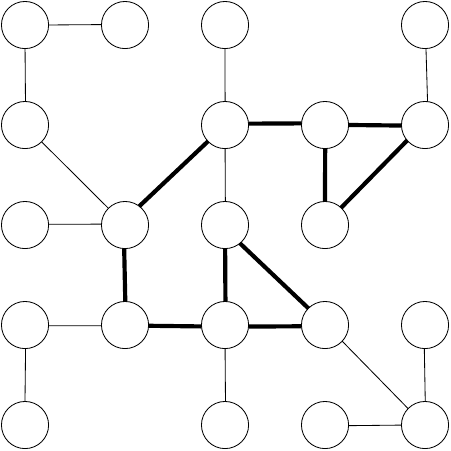}} & 
\parbox[c]{3.5cm}{\includegraphics[page=2,scale=0.6]{kernelexample.pdf}} & 
\parbox[c]{3cm}{\includegraphics[page=3, scale=0.6]{kernelexample.pdf}} \\
(a) & (b) & (c)
\end{tabular}
\end{center}

{\small \noindent Fig.\ 5. (a) A connected graph $G$ with 21 vertices and surplus 3 (the vertex labels are suppressed for clarity and the core edges are drawn in bold); (b) the core of $G$ (the core vertices of degree $3$ or more are drawn in bold); (c) the kernel of $G$}

\medskip

Fix a core $C$ with $\vset(C)\subset [n]$. 
There is a natural bijection between the set of connected graphs $G$ with $\vset(G)=[n]$ and $C(G)=C$, and the set of labelled forests $\forest_n^{\vset(C)}$ with vertex-set $[n]$ and root-set $\vset(C)$. For, given such a graph $G$, we form a forest in $\forest_n^{\vset(C)}$ by removing all edges of $C$ and rooting each connected component of the resulting graph at its unique element of $\vset(C)$. The bijection between trees and sequences described in Section~5 extends straightforwardly to forests, and from this we may deduce that if $|\vset(C)|=k$, then 
\[
|\{\mathrm{Graphs}~G~:~\vset(G)=[n], C(G)=C\}| = k n^{n-k-1}.
\]
From this, fairly standard combinatorial considerations yield the following result (see Proposition 1.5.1 of \cite{PIMSCRM} for a proof).

\medskip

\noindent \textbf{Theorem 7.1} \
\emph{Fix a kernel $K$ with vertex-set $[k]$ and with $m$ edges. 
Then as $\ell \to \infty$, the number of
cores $C$ with $\vset(C)=[\ell]$ and $K(C)=K$ is}
\[
(1+o(1))
\prod_{e \in \eset(K)} \frac{ 2^{-\mult(e;K)\I{\text{$e$ is a loop}}} }{\mult(e;K)!}\cdot \frac{(\ell-k)! \ \ell^{m-1}}{(m-1)!}.
\]

\medskip

From this, it follows that the kernel of a large random core with fixed surplus $s$ $(\ge 2)$ is whp $3$-regular. More specifically, fix $s \ge 2$ and let $C=C_\ell$ be a uniformly random core with $\vset(C)=[\ell]$ and surplus $s(C)=s$. Then, as $\ell \to \infty$,  $K(C)$ is whp a $3$-regular multigraph with $2(s-1)$ vertices and $3(s-1)$ edges, and for any fixed such multigraph $K$, 
 \begin{equation}\label{eq:kernel_dist}
\lim_{\ell \to \infty} \pr(K(C_\ell)=K) \propto \prod_{e \in \eset(K)} \frac{2^{-\I{\text{$e$ is a loop}}}}{\mult(e;K)!}.
\end{equation}
This in turn allows us to deduce that, for fixed $s \ge 2$, as $\ell \to \infty$, 
\begin{align*}
& |\{\mathrm{Cores}~C: \vset(C)=[\ell],s(C)=s\}|\\
& = (1+o(1))\frac{\ell! \ \ell^{3s-4}}{(2s-2)!(3s-4)!} \sum_{\substack{\mathrm{kernels}~K:\\ \vset(K)=[2(s-1)],\\ s(K)=s}}
\prod_{e \in \eset(K)} \frac{2^{-\I{\text{$e$ is a loop}}}} {\mult(e;K)!}.
\end{align*}
The case $s=1$ is much simpler, since a core of surplus $1$ is just a cycle: 
\[
|\{\mathrm{Cores}~C: \vset(C)=[\ell],s(C)=1\}|
= (\ell-1)!.
\]
This agrees in order of magnitude with the previous asymptotic value $\ell!\ \ell^{3s-4}$, on taking $s=1$. Combining this with the formula for the number of graphs with a given core, we may deduce the asymptotic number of connected graphs with given size and surplus, as well as the scaling of the core size in a large random connected graph with fixed surplus. In the next theorem, which is due to Wright~\cite{MR0463026} and Spencer~\cite{MR1431811}, we write
\[
\kappa(s) = \frac{1}{(2s-2)!(3s-4)!} \sum_{\substack{\mathrm{kernels}~K:\\ \vset(K)=[2(s-1)],\\ s(K)=s}}
\prod_{e \in \eset(K)} \frac{2^{-\I{\text{$e$ is a loop}}}}{\mult(e;K)!},
\]
for each integer $s \ge 2$, and set $\kappa(1)=1$. 

\medskip

\noindent \textbf{Theorem 7.2} \ 
\emph{For any fixed positive integer $s (\ge 1)$, as $n \to \infty$, 
\[
|\graphs_{n}^s| = (1+o(1)) \kappa(s) \cdot n^{n-2+3s/2} \int_0^\infty x^{3s-3}e^{-x^2/2}dx.
\]
Moreover, if $G \in_u \graphs_{n}^s$, then $|\vset(C(G))|/n^{1/2} \convdist X$, where the random variable $X$ has density 
\[
\frac{1}{2^{(3s-4)/2}\Gamma(\tfrac{1}{2}(3s-2))} x^{3s-3}e^{-x^2/2}\I{x \ge 0}.
\]
Here, $\Gamma$ denotes the gamma function. Equivalently, $X^2$ is distributed as \\$\mathrm{Gamma}(\tfrac{1}{2}(3s-2), \tfrac{1}{2})$.}

This theorem, together with the fact that the lengths of the paths in the core which are contracted to form kernel edges are exchangeable, has the remarkable consequence that, when $s \ge 2$, the asymptotic core structure has exactly the right distribution from which to `launch' a line-breaking construction. More precisely, a random metric space distributed as the scaling limit of a large random connected graph with surplus $s$ may be constructed as follows (see \cite{MR2650781}). 

\medskip

\noindent
\fbox{
\mbox{\begin{minipage}[t]{0.97\textwidth}
\textsc{Construction of a continuum random graph with fixed surplus $s$ $(\ge 2)$.}
\begin{enumerate}[itemindent=0pt,itemsep=0pt,leftmargin=15pt,topsep=2pt]
\item Let $K$ be a random $3$-regular multigraph sampled according to the distribution in \eqref{eq:kernel_dist}. Then sample $X\eqdist \smash{\sqrt{\mathrm{Gamma}\left(\tfrac{1}{2}(3s-2),\tfrac{1}{2}\right)}}$ and $(Y_1, Y_2,\ldots,Y_{3s-3})\eqdist\mathrm{Dirichlet}(1,1,\ldots,1)$, independently of each other and of $K$. Fix a labelling of the edges of $K$ by the elements of $[3s-3]$, and for $i \in [3s-3]$ replace the $i$th edge of $K$ by a line segment of length $X\cdot Y_i$.
\item Set $s_0=X$, and let $(s_i,i \ge 1)$ be the arrival times of an inhomogeneous Poisson process on the time-interval $[X,\infty)$ with rate $t$ at time $t$. Then, for each $i \ge 1$, attach a line segment of length $s_i-s_{i-1}$ to a point distributed according to the normalized Lebesgue measure on the object constructed so far. Finally, take the completion of the $i \to \infty$ limit to form a random compact metric space $(G^s, d^s)$. 
\end{enumerate}
\end{minipage}}
}

\medskip

A nearly identical construction works when $s=1$, but the line-breaking construction must be started with a `first line already attached' in order for the necessary distributional identities to work out.

\section{Constructing a graph from the depth-first queue process}
\setcounter{section}{8}
\setcounter{equation}{0}

\noindent In this section, we use the DFQ perspective to give a different (and complementary) way of obtaining the scaling limit of a uniformly random element of $\graphs_n^s$.  

We begin by arguing that we may associate a DFQ process with any connected graph $G \in \graphs_n$.  By convention, we take the vertex labelled 1 to be the root of the graph. We then explore $G$ in a depth-first fashion (just as we would if $G$ were a tree).  The `children' of the root are all of its neighbours, ordered by increasing vertex-label. We place these on a stack (as before) with the smallest-labelled vertex at the top. At each subsequent step, we explore (that is, reveal the neighbourhood of) whichever vertex sits at the top of the stack.  Its `children' are defined to be those of its neighbours in the graph which have not already been explored and are not already on the stack, again ordered by increasing vertex-label. We repeat this process until we have exhausted all of the vertices (and so the stack is empty). Along the way, we pick out a particular spanning tree $T(G)$ of the connected graph $G$, which we call the \emph{depth-first tree} of $G$.  The DFQ process of $G$ is then simply the DFQ process of $T(G)$.

Suppose now that we fix a particular labelled tree $T \in \graphs_n^0$, and ask which graphs $G$ have $T(G) = T$; these are necessarily graphs with $T$ as a spanning tree, along with a number of surplus edges. In order for these surplus edges to have been ignored by the DFQ process (and so for us to have obtained $T$ as our depth-first tree), they must have gone from the vertex $v$ explored at some particular step to another vertex which is on the stack at that time.  If $(q_i)_{0 \le i \le n}$ is the DFQ process of $T$, then there are, in total, $a(T) = \sum_{i=1}^{n-1}(q_i - 1)$ possible surplus edges that we may insert without altering the depth-first tree. (We call $a(T)$ the \emph{area} of the tree, since it is the area under the discrete function $(q_i)_{0 \le i \le n}$.) So there are $2^{a(T)}$ possible graphs $G$ with $T(G) = T$.  Since we have an ordered stack at each step, we can go further and identify each such graph $G$ with its DFQ process along with a subset of 
\[
S = \{(i,j): 0 \le i \le n, 1 \le j < q_i\}.
\]
Indeed, the DFQ process encodes the tree $T(G)$, and then a point at $(i,j)$ indicates that there is a surplus edge from the vertex explored at step $i$ to the vertex in position $j$ on the stack (counting from the bottom) (see Fig.\ 5).

\medskip 

\begin{center}
\begin{tabular}{cc}
\includegraphics[scale=0.58,page=1]{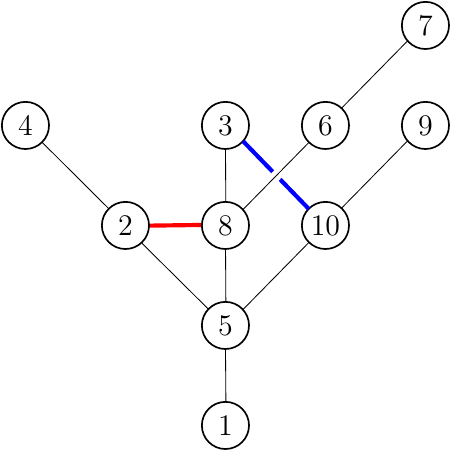} &
\includegraphics[scale=0.58,page=2]{DFQprocessmarks}\\
(a) & (b) \\
\end{tabular}
\end{center}
{\small \noindent Fig. 5. (a) A labelled graph, with its depth-first tree having lighter edges and its two surplus edges with heavy lines; (b) the corresponding DFQ process, with the possible locations for marks indicated by squares which are empty if the corresponding edge is absent, and filled if it is present}

\medskip

Now let us fix $s$ $(\ge 1)$ and see how this helps us to sample a uniform element of $\graphs_n^s$.  In order to do this, we need to pick uniformly from the collection of possible DFQ processes for $n$-vertex graphs which are endowed with sets of $s$ distinct points.  Any particular DFQ process $(q_i)_{0 \le i \le n}$ appears $\binom{a(T)}{s}$ times, once for each of the possible $s$-point sets with which we might decorate it. So equivalently, we may first sample a tree $T$ with probability proportional to $\binom{a(T)}{s}$, and then pick a uniform set of $s$ distinct points from $S$.  This procedure turns out to pass nicely to the limit. Indeed, the set of trees $T$ in which we are interested are those labelled by $[n]$, and a uniform choice from this set is equivalent to a size-conditioned Bienaym\'e tree with the Poisson distribution of mean 1 as its offspring distribution.  We do not want a uniform choice, but rather one biased by 
\[
\binom{a(T)}{s} = \frac{a(T) (a(T)-1) \ldots (a(T)-s+1)}{s!},
\]
and it turns out that we can control this bias. The uniform labelled tree lies in the Brownian domain of attraction, and so (in particular)  the DFQ process of such a random tree satisfies
\[
\frac{1}{\sqrt{n}} (q_{\lfloor nt \rfloor}, 0 \le t \le 1) \convdist 2\mathbf{e}, \quad \text{as $n \to \infty$.}
\]

Now $a(T) = \sum_{i=1}^n(q_{i} - 1) = n \int_0^1 (q_{\lfloor nt \rfloor} - 1) dt$, and so, by the continuous mapping theorem,
\[
\frac{1}{n^{3/2}} a(T) = \int_0^1 \frac{1}{\sqrt{n}}(q_{\lfloor nt \rfloor} - 1)\, dt \convdist 2 \int_0^1 \mathbf{e}(t) dt.
\]
Since $s$ is fixed, it follows (on checking the appropriate integrability) that if $\tilde{q}^s$ is the DFQ process of the biased tree, then
\[
\frac{1}{\sqrt{n}} (\tilde{q}^s_{\lfloor nt \rfloor}, 0 \le t \le 1) \convdist 2\tilde{\mathbf{e}}^s,
\]
where $\tilde{\mathbf{e}}^s$ is a continuous excursion whose distribution is determined by the fact that, for any bounded continuous functional $\Phi : C([0,1],\R^+) \to \R$,
\[
\E \left[\Phi(\tilde{\mathbf{e}}^s) \right] = \frac{\E\left[\Phi(\mathbf{e}) \left( \int_0^1 \mathbf{e}(t) dt \right)^s\right]}{ \E \left[\left( \int_0^1 \mathbf{e}(t) dt \right)^s \right]}.
\]
In other words, this is a Brownian excursion which has been biased by the $s$th power of its area.  So the scaling limit of the depth-first tree of the graph is the tree~$T_{2 \tilde{\mathbf{e}}}$.

In order to use this approach to derive the scaling limit of the graph itself, we also need to understand what happens to the points which tell us where to put the surplus edges. In the discrete picture, the points under the graph of the function $(\tilde{q}^s_i, 0 \le i \le n)$ and lying strictly above the $x$-axis correspond to possible surplus edges. We observe that the vertices on the stack at a particular step are all younger siblings of ancestors of the vertex being explored at that step -- in other words, they are at distance 1 from an ancestor of the vertex currently being explored, a distance that will be negligible once we rescale distances by $1/\sqrt{n}$.  So we may think of surplus edges as effectively going from the currently explored vertex to one of its ancestors. 

If we scale time by $1/n$ and space by $1/\sqrt{n}$, we then get convergence of our uniform $s$-set of distinct points in the discrete area to a collection of i.i.d.\ uniform points in the region under the function $2 \tilde{\mathbf{e}}^s$ and above the $x$-axis. As argued above, the DFQ process of a tree in the Brownian domain of attraction is asymptotically the same (up to a constant) as the height of its contour process.  This property is inherited when we bias the distribution, and so we obtain the result that each point $\{(x,y): 0\le x \le 1, 0\le y \le 2\tilde{\mathbf{e}}^s(x)\}$ represents a surplus `edge' of zero length, connecting the point at distance $2 \tilde{\mathbf{e}}^s(x)$ away from the root of $T_{2 \tilde{\mathbf{e}}^s}$, which corresponds to $x$, and a second point that is at distance $y$ away from the root along the same ancestral line. Formally, we identify each of the $s$ pairs of points in order to obtain a limiting metric space $(G^s, d^s)$, which we endow with the probability measure $\mu^s$ that is the push-forward of the measure $\mu_{2 \tilde{\mathbf{e}}^s}$ under the identifications. That this gives the same object as the one described in the previous section is a far-from-obvious extension of the equivalence of the two definitions of the Brownian CRT, one from the coding excursion and the other via line-breaking.

The following result is due to Addario-Berry, Broutin and Goldschmidt~\cite{MR2892951}.

\medskip

\noindent \textbf{Theorem 8.1} \
\emph{Fix $s \ge 1$, and let $G_n^s$ be a uniformly random element of $\graphs_n^s$. Let $d_n^s$ be the graph distance on $G_n^s$, and $\mu_n^s$ be the uniform measure on the vertices. Then
\[
(V(G^s_n), n^{-1/2} d_n^s, \mu_n^s) \convdist (G^s, d^s, \mu^s),
\]
as $n \to \infty$ in the GHP sense.}

\medskip

We conclude this section with an observation which is useful in the next one. Suppose that we wish to understand the relationship between the scaling limit of an element of $\graphs_n^s$ and an element of $\graphs_{\lfloor nx\rfloor}^s$. The former has a CRT encoded by $\tilde{\mathbf{e}}^s$ as the scaling limit of its depth-first tree. For the latter, we may argue in the same way as before to see that its depth-first tree has scaling limit encoded by an excursion $\tilde{\mathbf{e}}^s_x$ defined by
\[
\E \left[\Phi(\tilde{\mathbf{e}}^s_x) \right] = \frac{\E\left[\Phi(\mathbf{e}_x) \left( \int_0^x \mathbf{e}_x(t) dt \right)^s\right]}{ \E \left[\left( \int_0^x \mathbf{e}_x(t) dt \right)^s \right]},
\]
for any bounded continuous functional $\Phi: C([0,x], \R^+) \to \R$. By Brownian scaling (see \eqref{eq:bscale}), this is equivalent to the expression
\[
\E \left[\Phi(\tilde{\mathbf{e}}^s_x) \right] = \frac{\E\left[\Phi(x^{1/2}\mathbf{e}(x \cdot)) \left(\int_0^1 \mathbf{e}(t) dt \right)^s\right]}{ \E \left[\left( \int_0^1 \mathbf{e}(t) dt \right)^s \right]} = \E\left[\Phi(x^{1/2} \tilde{\mathbf{e}}^s(\cdot/x))\right].
\]
In other words, we get a time- and space-rescaled copy of the unit length biased excursion, which encodes the tree $T_{2\tilde{\mathbf{e}}^s}$ with all its distances multiplied by $\sqrt{x}$.

\section{The critical Erd\H{o}s--R\'enyi random graph}

\setcounter{section}{9}
\setcounter{equation}{0}

\noindent An application of the results of the previous section may be found in the context of the Erd\H{o}s--R\'enyi random graph $\mathbb{G}_{n,p}$. We begin with a simple observation: given its vertex-set and number of surplus edges, any component of an Erd\H{o}s--R\'enyi random graph is a uniformly random labelled graph on that vertex-set with the given number of surplus edges. Moreover, the different components are conditionally independent, given their vertex-sets.  In the critical setting, when $p = 1/n + \lambda n^{-4/3}$, for some $\lambda \in \R$, it turns out that the largest components have sizes of the order of $n^{2/3}$ whp, and crucially the surpluses of those components are $\Theta(1)$. So if we can understand the distributions of those sizes and surpluses, we may then simply apply the results of the previous section in order to obtain the metric space scaling limit of the whole graph. The result on sizes and surpluses that we need is due to Aldous~\cite{MR1434128}. (This result was one of the first instances of the application of the tools of stochastic analysis in the context of random graphs.)

Fix $p = 1/n + \lambda n^{-4/3}$, and let $C_{n,1}, C_{n,2}, \ldots$ be the components of $\mathbb{G}_{n,p}$, listed in decreasing order of size, and let $S_{n,1}, S_{n,2}, \ldots$ be the surpluses of those components. Let $\ell_2^{\downarrow} = \{(x_1, x_2, \ldots): x_1 \ge x_2 \ge \ldots \ge 0, \sum_{i=1}^{\infty} x_i^2 < \infty\}$. Now let $(B_t)_{t \ge 0}$ be a standard Brownian motion, and define
\[
B^{\lambda}_t = B_t + \lambda t - t^2/2, \quad t \ge 0.
\]
We reflect this process at its running infimum to obtain a process $(R_t)_{t \ge 0}$ which stays non-negative for all time -- that is, $R_t = B^{\lambda}_t - \inf_{0 \le s \le t} B^{\lambda}_s$.  Now, conditionally on $R$, let $(P_t)_{t \ge 0}$ be a Poisson process of intensity $R_t$ at time $t$. Finally, let $\gamma_1, \gamma_2, \ldots$ be the lengths of the excursions of $R$ above 0, listed in decreasing order, and let $\sigma_1, \sigma_2, \ldots$ be the number of arrivals of the Poisson process in each of those excursions. (It may be shown that $\sum_{i=1}^{\infty} \gamma_i^2 < \infty$ almost surely.) Then Aldous's result is as follows.

\medskip
\noindent \textbf{Theorem 9.1} \ 
\emph{As $n \to \infty$,
\[
n^{-2/3}(|V(C_{n,1})|, |V(C_{n,2})|, \ldots) \convdist (\gamma_1, \gamma_2, \ldots)
\]
in $\ell_2^{\downarrow}$, endowed with the $\ell_2$-norm, jointly with the convergence
\[
(S_{n,1}, S_{n,2}, \ldots) \convdist (\sigma_1, \sigma_2, \ldots),
\]
in the sense of the product topology.}

\medskip

Let $d_{n,i}$ be the graph distance on $C_{n,i}$, and let $\mu_{n,i}$ be the uniform measure on the vertices. In view of the above discussion, we now obtain a scaling limit theorem for the critical Erd\H{o}s--R\'enyi random graph (see Addario-Berry, Broutin and Goldschmidt~\cite{MR2892951}).

\medskip
\noindent \textbf{Theorem 9.2} \
\emph{Jointly with the convergence in Theorem 9.1, we have
\[
((V(C_{n,i}), n^{-1/3} d_{n,i}, \mu_{n,i}), i \ge 1) \convdist ((G_i^{\sigma_i}, \sqrt{\gamma_i} d_i^{\sigma_i}, \mu_i^{\sigma_i}), i \ge 1),
\]
as $n \to \infty$, in an $\ell_4$-sequence version of GHP, where, for each $s \ge 0$, the entries of $(G_i^s, d_i^s, \mu_i^s), i \ge 1$ are i.i.d.\ copies of $(G^s, d^s ,\mu^s)$, independent of $(\gamma_1, \gamma_2, \ldots)$ and $(\sigma_1, \sigma_2, \ldots)$.}

\medskip

Theorem 9.1 is perhaps surprising at first sight, so let us try to give some insight into why it is true; the key to understanding it is again the DFQ process. We explore the components of the graph one by one, starting from the vertex labelled 1. We find the DFQ process associated with the first component, and then successively concatenate onto its end the DFQ processes associated with the other components, each time choosing which one to explore by starting from the lowest-labelled non-explored vertex in the whole graph.  

This concatenated process $(Q_i)_{0 \le i \le n}$ is actually Markovian with relatively simple transition probabilities: if $Q_i = 0$ then $Q_{i+1} = 1$, since we start a new component; but if $Q_i = q > 0$, then $Q_{i+1} - Q_i \eqdist \Bin(n- i - q , 1/n + \lambda n^{-4/3}) - 1$. In order to see this, note that at step $i$ we have already explored $i$ vertices and have $q$ more on the stack. So there are $n-i-q$ other unexplored vertices in the graph, and there is an edge between the vertex $v_{i+1}$ explored at step $i+1$ and each of these vertices independently with probability $p = 1/n + \lambda n^{-4/3}$. So $v_{i+1}$ has a $\Bin(n- i - q , 1/n + \lambda n^{-4/3})$ number of children, which yields the claimed transition probability for the DFQ process. 

We then have
\begin{align*}
\E\left[Q_{i+1} - Q_i | Q_i = q \right] & = \left(\frac{1}{n} + \frac{\lambda}{n^{4/3}} \right) (n - i - q) - 1 \\
& =  \frac{\lambda}{n^{1/3}} - \frac{i}{n} + o(n^{-1/3}),
\end{align*}
as long as $q > 0$, $q = o(n^{2/3})$, and $i = o(n)$. In other words, when it is away from 0, $Q$ is essentially a random walk with a slight time-dependent drift, which should therefore look asymptotically like a Brownian motion with drift.  Taking $i = \lfloor n^{2/3} t \rfloor$, we see that
\[
\E\left[Q_{\lfloor n^{2/3} t \rfloor + 1} - Q_{\lfloor n^{2/3} t \rfloor} \big| Q_{\lfloor n^{2/3} t \rfloor} = q \right] = \frac{1}{n^{1/3}} (\lambda - t) + o(n^{-1/3})\, ,
\]
whenever $q=o(n^{2/3})$. 
Recalling that the process reflects at 0, this line of argument eventually yields the result that
\[
n^{-1/3}(Q_{\lfloor n^{2/3} t \rfloor}, t \ge 0) \convdist R.
\]
Since the sizes of the components are encoded in $Q$ as the lengths of excursions away from 0, we see that these are indeed $O(n^{2/3})$, and have the lengths of the excursions of $R$ away from 0 as their scaling limits.  In order to get the numbers of surplus edges, it suffices to recall that surplus edges arise as edges joining the vertex currently being explored to a vertex on the stack. It follows that, at step $i$, and conditionally on the stack having size $Q_i = q$, there are $\Bin(q-1, 1/n + \lambda n^{-4/3})$ surplus edges, with expectation given approximately by $q/n = n^{-2/3} (q/n^{1/3})$. Our rescaling is now exactly right for us to obtain the claimed Poisson process of surplus edges in the limit. 

\section{Uniform random graphs with i.i.d.\ degrees}

\setcounter{section}{10}
\setcounter{equation}{0}

\noindent Just as the Brownian CRT is the universal scaling limit of a broad family of different random tree models, so the objects in the previous two sections arise outside the context of the Erd\H{o}s--R\'enyi random graph.  There are two particular models which have been particularly well studied from the perspective of scaling limits at criticality: the uniform random graph with given degrees, and the muliplicative random graph (which goes by various different names, including the `rank-1 inhomogeneous random graph' and the `Poisson graph').  In the interests of brevity, we focus our attention on the uniform random graph with given i.i.d.\ random degrees, which was briefly described in Section~2.  

More precisely, let $D_1, D_2, \ldots, D_n$ be i.i.d.\ copies of a strictly positive integer-valued random variable $D$ for which $\pr(D = 2) < 1$ and $\E[D^2] < \infty$; these represent the degrees of vertices $1, 2, \ldots, n$, respectively. In general, there is no reason why $\sum_{i=1}^n D_i$ should be even; if it is odd, we simply reduce the degree of vertex $n$ by 1. With this modification, it turns out that there is an asymptotically strictly positive probability that a simple graph with these degrees exists so that, given $D_1, D_2, \ldots, D_n$ and taking $n$ large enough, it makes sense to sample a uniformly random graph with these degrees. Define $\mu = \E[D]$ and $\theta = \E[D(D-1)]/\E[D]$. Then if $\theta < 1$, the largest component contains a vanishing proportion of the vertices, whereas if $\theta > 1$, the largest component contains a strictly positive proportion of the vertices, whp as $n \to \infty$.  So a critical graph is obtained by taking $\theta = 1$.  

Suppose now that $\E[D^3] < \infty$, and define $\beta = \E[D(D-1)(D-2)]$. Then the precise analogues of Theorems 9.1 and 9.2 hold for the components of this graph, with the small modification that we replace $(B^{\lambda}_t)_{t \ge 0}$ by $(\tilde{B}_t)_{t \ge 0}$, where
\[
\tilde{B}_t = \sqrt{\beta/\mu} \, B_t - (\beta/2 \mu^2) t^2.
\]

If the condition $\E[D^3] < \infty$ fails, but we have regular tail-behaviour of the form 
\[
\pr(D = k) = (1+o(1))c k^{-(\alpha+2)} \quad \text{as $k \to \infty$,}
\]
for some $c > 0$ and $\alpha \in (1,2)$, then the connected components of the resulting graph resemble $\alpha$-stable CRTs, rather than Brownian CRTs. Indeed, a complete analogue of the theory we have described for the Erd\H{o}s--R\'enyi random graph exists in that setting also, including both functional encodings and line-breaking constructions (see \cite{MR4515689}, \cite{MR4526241}).

\fontsize{10}{12}
\selectfont
\small
\bibliographystyle{plainnat}

\normalsize

\end{document}